\documentstyle[amstex,amscd,euscript,amssymb]{amsart}
\numberwithin{equation}{section}
\newtheorem{thm}{Theorem}[section]

\theoremstyle{remark}

\begin{document}
\newcommand{\dcp}{\,{\triangleright \kern-0.145em \triangleleft}}
\newcommand{\dcpP}{\,{\triangleright \kern-0.145em \triangleleft}_{\cal P}}
\newcommand{\tensE}{{\otimes}_{E}}
\newcommand{\tenc}{{\intercal}}
\newcommand{\cu}{\varepsilon}
\newcommand{\enal}{e_{\lambda}}
\newcommand{\enam}{e_{\mu}}
\newcommand{\enan}{e_{\nu}}
\newcommand{\enax}{e_{\xi}}
\newcommand{\ema}{{\stackrel{\scriptscriptstyle\circ}{e}}}
\newcommand{\emal}{{\stackrel{\scriptscriptstyle\circ}{e}_{\lambda}}}
\newcommand{\emam}{{\stackrel{\scriptscriptstyle\circ}{e}_{\mu}}}
\newcommand{\eman}{{\stackrel{\scriptscriptstyle\circ}{e}_{\nu}}}
\newcommand{\emax}{{\stackrel{\scriptscriptstyle\circ}{e}_{\xi}}}
\newcommand{\emai}{{\stackrel{\scriptscriptstyle\circ}{e}_i}}
\newcommand{\wsta}{{\stackrel{\scriptscriptstyle *}{w}}}
\newcommand{\ehoi}{e_{{\HH}^{\circ},i}}
\newcommand{\ehoj}{e_{{\HH}^{\circ},j}}
\newcommand{\ehok}{e_{{\HH}^{\circ},k}}
\newcommand{\ehol}{e_{{\HH}^{\circ},l}}
\newcommand{\emahoi}{{\stackrel{\scriptscriptstyle\circ}{e}_{{\HH}^{\circ},i}}}
\newcommand{\emahoj}{{\stackrel{\scriptscriptstyle\circ}{e}_{{\HH}^{\circ},j}}}
\newcommand{\emahok}{{\stackrel{\scriptscriptstyle\circ}{e}_{{\HH}^{\circ},k}}}
\newcommand{\emahol}{{\stackrel{\scriptscriptstyle\circ}{e}_{{\HH}^{\circ},l}}}
\newcommand{\epma}{{\stackrel{\scriptscriptstyle\circ}{\varepsilon}}}
\newcommand{\epmai}{{\stackrel{\scriptscriptstyle\circ}{\varepsilon}_i}}
\newcommand{\epmaj}{{\stackrel{\scriptscriptstyle\circ}{\varepsilon}_j}}
\newcommand{\epmak}{{\stackrel{\scriptscriptstyle\circ}{\varepsilon}_k}}
\newcommand{\epmal}{{\stackrel{\scriptscriptstyle\circ}{\varepsilon}_l}}
\newcommand{\lamma}{{\stackrel{\scriptscriptstyle\circ}{\lambda}}}
\newcommand{\muma}{{\stackrel{\scriptscriptstyle\circ}{\mu}}}
\newcommand{\numa}{{\stackrel{\scriptscriptstyle\circ}{\nu}}}
\newcommand{\xima}{{\stackrel{\scriptscriptstyle\circ}{\xi}}}
\newcommand{\Vma}{{\stackrel{\scriptscriptstyle\circ}{{\EuScript V}}}}
\newcommand{\Rma}{{\stackrel{\scriptscriptstyle\circ}{R}}}
\newcommand{\Hom}{\mathrm{Hom}}
\newcommand{\End}{\mathrm{End}}
\newcommand{\Tr}{\mathrm{Tr}}
\newcommand{\id}{\mathrm{id}}
\newcommand{\card}{\mathrm{card}}
\newcommand{\spa}{\mathrm{span}}
\newcommand{\op}{\mathrm{op}}
\newcommand{\cop}{\mathrm{cop}}
\newcommand{\bop}{\mathrm{bop}}
\newcommand{\ob}{\mathrm{ob}}
\newcommand{\bimod}{\bold{bmd}}
\newcommand{\bmdr}{{{\bold{bmd}}} ({\it R})}
\newcommand{\bmdrma}{{\bold{bmd}}
({{\stackrel{\scriptscriptstyle\circ}{{\it R}}}})}
\newcommand{\bmde}{\bold{Bmd} (E)}
\newcommand{\HH}{\frak H}
\newcommand{\KK}{\frak K}
\newcommand{\Hhat}{\hat{\frak H}}
\newcommand{\Ht}{\tilde{\frak H}}
\newcommand{\HGm}{\frak{H} [G^{-1}]}
\newcommand{\Hc}{\mathrm{Hc}}
\newcommand{\HcH}{\mathrm{Hc}({\frak H})}
\newcommand{\Aw}{{\frak A}(w)}
\newcommand{\SS}{{\frak S}(A_{N-1};t)_{\epsilon}}
\newcommand{\SSe}{{\frak S}(A_{N-1};t)_{\epsilon,\zeta}}
\newcommand{\SSz}{{\frak S}(A_{N-1};t)_{\epsilon,\zeta}}
\newcommand{\SSi}{{\frak S}(A_{N-1};t)_{\epsilon,\zeta}^{\iota}}
\newcommand{\Ss}{{\frak S}(A_{1};t)_{\epsilon}}
\newcommand{\AS}{{\frak A}(w_{N,t,\epsilon})}
\newcommand{\ASe}{{\frak A}(w_{N,t,\epsilon,\zeta})}
\newcommand{\Str}[1]{\mathrm{Str}^{#1} (\G, *)}
\newcommand{\Len}{\EuScript{L}}
\newcommand{\sgn}[1]{(- \epsilon)^{\EuScript{L} ( #1 )}}
\newcommand{\Rp}{{\cal R}^+}
\newcommand{\Rm}{{\cal R}^-}
\newcommand{\Rpm}{{\cal R}^{\pm}}
\newcommand{\Rmp}{{\cal R}^{\mp}}
\newcommand{\Rph}{{\hat{\cal R}}^+}
\newcommand{\Rmh}{{\hat{\cal R}}^-}
\newcommand{\Rpmh}{{\hat{\cal R}}^{\pm}}
\newcommand{\Rmph}{{\hat{\cal R}}^{\mp}}
\newcommand{\Rt}{{\tilde{\cal R}}}
\newcommand{\Rpt}{{\tilde{\cal R}}^+}
\newcommand{\Rmt}{{\tilde{\cal R}}^-}
\newcommand{\Rpmt}{{\tilde{\cal R}}^{\pm}}
\newcommand{\G}{{\EuScript G}}
\newcommand{\V}{{\EuScript V}}
\newcommand{\K}{{\Bbb K}}
\newcommand{\C}{{\Bbb C}}
\newcommand{\R}{{\Bbb R}}
\newcommand{\Q}{\Bbb Q}
\newcommand{\Z}{{\Bbb Z}}
\newcommand{\LD}{\mathrm{LD}}
\newcommand{\GLD}{\G_{\mathrm{LD}}}
\newcommand{\wLD}{w_{\mathrm{LD}}}
\newcommand{\KP}{\mathrm{KP}}
\newcommand{\GLEH}{\mathrm{GLE} (\HH)}
\newcommand{\LLam}{\boldsymbol\Lambda}
\newcommand{\aaa}{\bold a}
\newcommand{\bbb}{\bold b}
\newcommand{\ccc}{\bold c}
\newcommand{\ddd}{\bold d}
\newcommand{\p}{\bold p}
\newcommand{\q}{\bold q}
\newcommand{\r}{\bold r}
\newcommand{\s}{\bold s}
\newcommand{\taaa}{\tilde{\bold a}}
\newcommand{\tbbb}{\tilde{\bold b}}
\newcommand{\tccc}{\tilde{\bold c}}
\newcommand{\tddd}{\tilde{\bold d}}
\newcommand{\tp}{\tilde{\bold p}}
\newcommand{\tq}{\tilde{\bold q}}
\newcommand{\tr}{\tilde{\bold r}}
\newcommand{\ts}{\tilde{\bold s}}
\newcommand{\st}{\frak {s}}
\newcommand{\en}{\frak {r}}
\newcommand{\suma}{\sum_{(a)}}
\newcommand{\sumb}{\sum_{(b)}}
\newcommand{\sumc}{\sum_{(c)}}
\newcommand{\sumd}{\sum_{(d)}}
\newcommand{\sumx}{\sum_{(x)}}
\newcommand{\sumy}{\sum_{(y)}}
\newcommand{\sumab}{\sum_{(a),(b)}}
\newcommand{\sumbc}{\sum_{(b),(c)}}
\newcommand{\sumad}{\sum_{(a),(d)}}
\newcommand{\sumxy}{\sum_{(x),(y)}}
\newcommand{\sumax}{\sum_{(a),(x)}}
\newcommand{\sumbx}{\sum_{(b),(x)}}
\newcommand{\sumabc}{\sum_{(a),(b),(c)}}
\newcommand{\sumabcd}{\sum_{(a),(b),(c),(d)}}
\newcommand{\sumk}{\sum_{k \in \V}}
\newcommand{\suml}{\sum_{l \in \V}}
\newcommand{\sumkl}{\sum_{k,l \in \V}}
\newcommand{\face}[4]{\left( \scriptstyle{#1} 
                      \textstyle{\frac[0pt]{#2}{#3}} \scriptstyle{#4} \right)}
\newcommand{\Wrpqs}[4]{w \! \left[ #1 \, \frac[0pt]{#2}{#3} \, #4 \right]}
\newcommand{\wrpqs}[4]{w \! \left[ \scriptstyle{#1} 
                      \textstyle{\frac[0pt]{#2}{#3}} \scriptstyle{#4} \right]}
\newcommand{\hijk}[4]{ \!\! \left[ {#1 \,\, #2} \atop {#3 \,\, #4} \right]}

\title[A canonical Tannaka duality]{A canonical Tannaka duality for finite
semisimple tensor categories}
\author[Takahiro Hayashi]{Takahiro Hayashi \\
$\quad$\\
${\rm Graduate}\, {\rm School}\, {\rm of}\, {\rm Mathematics},\, 
{\rm Nagoya}\, {\rm University,}$\\ 
$
{\rm Chikusa} \frac{\,\,}{\,\,}{\rm ku,}\, 
{\rm Nagoya}\,\, {\rm464}\frac{\,\,}{\,\,}{\rm 8602,}\, 
{\rm Japan}$}
\maketitle

\begin{abstract}	
 For each finite semisimple tensor category, we associate 
 a quantum group (face algebra) whose comodule category is
 equivalent to the original one, in a simple natural manner.
 To do this, we also give a generalization of the Tannaka-Krein
 duality, which assigns a face algebra for each tensor category
 equipped with an embedding into a certain kind of bimodule category.
\end{abstract}

\section{introduction}

Let ${\bold C}$ be a tensor category over a field $\K$
and $\Omega$ an ``embedding" of 
${\bold C}$ into the category $\bold{mod}(\K)$ of all 
finite-dimensional vector spaces over $\K$.
A version of the Tannaka-Krein duality theorem says that
there exists a natural way to construct a bialgebra 
$C (\Omega)$ such that the category $\bold{com} (C (\Omega))$ of all 
finite-dimensional right $C (\Omega)$-comodules is
equivalent to ${\bold C}$.

After the discovery of the quantum groups and Jones' index theory,  
it is recognized that there exist natural interesting examples of 
tensor categories which have no embeddings into $\bold{mod}(\K)$
(see e.g. \cite{GelfandKazhdan}, \cite{Ocneanu}, \cite{Turaev}).

In this letter, we show that
there exists a natural simple way to construct an embedding
$\Omega_0$ of each tensor category ${\bold C}$ into
a certain bimodule category $\bmdr$,
provided that ${\bold C}$ is finite
semisimple. 
Here $R$ denotes a suitable finite direct product of copies of $\K$.
Also, we give a generalization of the Tannaka-Krein duality 
to tensor categories ${\bold C}$ equipped with such embeddings 
$\Omega\!: {\bold C} \to \bmdr$, using a generalization of the 
notion of a bialgebra, which is called a {\it face algebra}
(cf. \cite{subf}-\cite{fb}, and also \cite{BS}, 
\cite{JurcoSchupp}, \cite{Schauenburgface}). 
Combining these, we obtain a natural way
to construct a face algebra $C(\Omega_0)$ such that
$\bold{com} (C (\Omega_0))$ is equivalent to ${\bold C}$,
where ${\bold C}$ is an arbitrary finite semisimple 
tensor category over $\K$. 

In Section 2, we recall the definition of the $\V$-face algebra
and that of the $\V$-{\it dressed coalgebra} (cf. \cite{fb}).
The latter enables us to imitate Schauenburg's formulation of 
the Tannaka-Krein duality
(cf. \cite{SchauenburgTD}).
In Section 3, we give a generalization of the Tannaka-Krein 
duality, which we call $\V$-{\it Tannaka duality}.
In Section 4, we state the main result.

Throughout this letter, we denote by 
$\bold{Mod} (A)$ and $\bold{mod}  (A)$
the category of all left $A$-modules 
and the category of all finite-dimensional $A$-modules
of a $\K$-algebra $A$, respectively, 
where $\K$ denotes the ground field. 
In particular, $\bold{Mod} (\K)$ denotes
the category of all vector spaces over $\K$.  
Also we denote by $\bold{bmd} (A)$ 
and $\bold{bmd} (A)$ the category of all $A$-bimodules 
and the category of all $\K$-finite-dimensional $A$-bimodules, 
respectively.
For a coalgebra $C$ and its right comodule $M$, we use Sweedler's
sigma notation, such as 
$\Delta (c)$ $=$ $\sum_{(c)} c_{(1)} \otimes c_{(2)}$
$(c \in C)$ and 
$(\mathrm{id} \otimes \Delta) \circ \delta (m)$
$=$ $\sum_{(m)} m_{(0)} \otimes m_{(1)} \otimes m_{(2)}$
$(m \in M)$, where $\Delta$ and $\delta$ denotes the coproduct of $C$
and the coaction of $M$, respectively (cf \cite{Sweedler}).
We denote by $\bold{com} (C)$ the category of all finite-dimensional
right $C$-comodules.

We wish to thank Prof. S. Yamagami for valuable discussions.

\section{Preliminary}

%
%
Let $\K$ be a field and ${\EuScript V}$ a finite non-empty set.
We denote  by $R = R_{\V}$ the $\K$-linear span of 
$\V$ equipped with an algebra structure given by
$\lambda \mu = \delta_{\lambda \mu} \lambda$
$(\lambda, \mu \in \V)$.
Let $\Rma = \bigoplus_{\lambda \in \V} \K \lamma$ be a copy of $R$ and 
$E_{\V} = E  = \Rma \otimes R$ the tensor product algebra equipped with
a coalgebra structure given by setting
\begin{equation}			
 \Delta(\lamma \mu) = \sum_{\nu \in \V} \lamma \nu \otimes \numa \mu,
\quad
 \cu (\lamma \mu) = \delta_{\lambda \mu}
\label{D(ee)}
\end{equation}
for each $\lambda, \mu \in \V$. 
Here we write $\lambda$ and $\lamma$ instead of 
$\mathrm{id} \otimes \lambda$ and $\lamma \otimes \mathrm{id}$,
respectively.
Let $\HH$ be a ${\Bbb K}$-algebra equipped
with a coalgebra structure $(\HH,\Delta,\cu)$ together with an
algebra-coalgebra map $E \to \HH$.
We say that $\HH$ is a {\em ${\V}$-face algebra} 
if the following relations are satisfied:
\begin{equation}
 \Delta (ab) = \Delta(a) \Delta(b),
\label{D(ab)}
\end{equation}
%
%
\begin{equation}
 \cu (ab) = \sum_{\nu \in \V} \cu (a\nu) \cu (\numa b)
\label{e(ab)}
\end{equation}
for each $a, b \in \HH$.
Here we denote the images of $\nu$ and $\numa$ of $E$ via 
the map $E \to \HH$ simply by $\nu$ and $\numa$, respectively.
We call the elements $\nu, \numa \in \HH$ {\it face idempotents}
of $\HH$. 
It is known that a bialgebra is an equivalent notion of 
a $\V$-face algebra with $\mathrm{card}(\V) =$ $1$.
%


We say that a linear map $S\!: {\HH} \to {\HH}$ is an {\em antipode} of ${\HH}$, 
or $({\HH}, S)$ is a {\em Hopf} ${\V}$-{\em face algebra} if 
\begin{equation}		
 \suma S(a_{(1)})a_{(2)} = \sum_{\nu \in \V} \cu (a \nu) \nu, 
 \quad				
 \suma a_{(1)}S(a_{(2)}) = \sum_{\nu \in \V} \cu (\nu a){\numa},
 \label{S(a)a}
\end{equation}
\begin{equation}
 \suma S(a_{(1)})a_{(2)}S(a_{(3)}) = S (a) 
 \label{S(a)aS(a)}
\end{equation} 
for each $a \in \HH$.
An antipode of a $\V$-face algebra is 
an antialgebra-anticoalgebra map, which satisfies
\begin{equation}
 S(\lamma \mu) = \muma \lambda 
 \quad (\lambda, \mu \in {\V}).
\label{S(ee)}
\end{equation}
%
%
%
The antipode of a $\V$-face algebra is unique 
if it exists.

%
%

Let $M$ be a right comodule of a $\V$-face algebra $\HH$.
We define an $R$-bimodule structure on $M$ by
\begin{equation}
 \lambda m \mu	=		
 \sum_{(m)} m_{(0)} \cu (\lambda m_{(1)} \mu)
 \quad (m \in M,\, \lambda, \mu \in \V).
\label{}
\end{equation}
Let $N$ be another $\HH$-comodule. We define an $\HH$-comodule structure
on $M \otimes_{R} N$ by
\begin{equation}
 m \otimes_R n	\mapsto		
 \sum_{(m), (n)} ( m_{(0)} 
 \otimes_R n_{(0)} ) \otimes
 m_{(1)} n_{(1)}
 \quad (m \in M,\, n \in N). 
\label{}
\end{equation}
The category $\bold{com} (\HH)$ 
becomes a monoidal category via this operation. 
If, in addition, $\HH$ has a bijective antipode, then  
$\bold{com} (\HH)$ is rigid and the left dual object of 
$M \in \ob \bold{com} (\HH)$ is given by
\begin{gather}
 M^{\vee} = \Hom_{\K} (M, \K), \nonumber\\
 \sum_{(m^{\vee})}
 \langle m^{\vee}_{(0)},\, m \rangle m^{\vee}_{(1)}
 =
 \sum_{(m)}
 \langle m^{\vee},\, m_{(0)} \rangle S(m_{(1)}) \nonumber\\
 (m \in M,\, m^{\vee} \in M^{\vee}). 
\label{Mveedef}
\end{gather}

Let $C$ be a coalgebra equipped with $E$-bimodule structure.
We say that $C$ is a $\V$-{\it dressed coalgebra} if $C$ satisfies 
%
\begin{equation} 			
\label{D(eeaee)} 
 \Delta (\lamma \mu c \lamma{}^{\prime} \mu{}^{\prime}) =
 \sumc \lamma c_{(1)} \lamma{}^{\prime} \otimes \mu c_{(2)} \mu{}^{\prime}, 
\end{equation}
\begin{equation}
 \sumc \lambda  c_{(1)} \mu \otimes c_{(2)} = 
 \sumc c_{(1)} \otimes \lamma c_{(2)} \muma,  
\label{eae*a}
\end{equation}
\begin{equation}			
 \cu (\lamma c \muma) = \cu (\lambda c \mu)
\label{e(ee)}
\end{equation}
for each $c \in C$ and $\lambda, \mu, \lambda', \mu' \in \V$.
For a $\V$-dressed coalgebra $C$, we define another 
$\V$-dressed coalgebra $C^{\bop}$ to be the opposite 
coalgebra of $C$ equipped with $E$-bimodule structure given by
$\lamma \mu \otimes c \otimes \lamma {}^{\prime} \mu^{\prime}$ $\mapsto$
$\muma{}^{\prime} {\lambda}^{\prime}  c \muma \lambda$
$(\lambda, \mu, \lambda {}^{\prime}, \mu^{\prime} \in \V)$.  
Let $D$ be another $\V$-dressed coalgebra. 
A coalgebra map $f$ from $C$ to $D$ is called a 
{\it map of $\V$-dressed coalgebras} if it is an
$E$-bimodule map.  
The tensor product $E$-bimodule
$C \tensE D$ becomes a $\V$-dressed coalgebra via 
\begin{equation} 			
\label{D(ctensEd)}
 \Delta (c \tensE d) =
 \sum_{(c), (d)} (c_{(1)} \tensE d_{(1)}) \otimes 
 (c_{(2)} \tensE d_{(2)}),
\end{equation}
\begin{equation} 			
\label{cu(ctensEd)}
 \cu (c \tensE d) =
 \sum_{\nu \in \V} \cu (c \nu) \cu (\nu d).
 \quad 
\end{equation}
The category $\bold{DrCoalg}_{\V}$ of all $\V$-dressed coalgebras
becomes a monoidal category
via this operation  with unit object $E$. 
It is easy to see that a monoid in this category is simply a $\V$-face algebra. 
%

%
%
%
%
%
%
%

%
%

\section{The $\V$-Tannaka duality}

For each $M \in \ob \bold{bmd} ({\it R})$ and $N \in \ob \bmdrma$,
we regard $N \otimes M$ as an 
object of $\bmde$ via  
\begin{equation}
 \lamma \mu (n \otimes m) \lamma{}^{\prime} \mu{}^{\prime} =  
 \lamma n \lamma {}^{\prime} \otimes \mu m \mu{}^{\prime}
 \quad (n \in N, m \in M, \lambda, \mu, \lambda {}^{\prime}, \mu{}^{\prime} \in \V).
\end{equation}
For each $M \in \mathrm{ob} \bmdr$, we define
$M^*$ to be the $\K$-linear dual of $M$ equipped with
$\Rma$-bimodule structure given by
\begin{equation}
 \langle \lamma m^* \muma,\, m \rangle =
 \langle m^*,\, \lambda m \mu \rangle 
 \quad (m \in M, m^* \in M^*, \lambda, \mu \in \V). 
\end{equation}
Also, we define
$M^{\vee}$ to be the $\K$-linear dual of $M$ equipped with
$R$-bimodule structure given by
\begin{equation}
 \langle \lambda m^{\vee} \mu,\, m \rangle =
 \langle m^{\vee},\, \mu m \lambda \rangle 
 \quad (m^{\vee} \in M^{\vee},\, m \in M,\, \lambda, \mu \in \V). 
\end{equation}

We say that a pair $({\bold C}, \Omega)$ is a {\it category with
$\V$-face} if ${\bold C}$ is an essentially small category and 
$\Omega$ is a functor from ${\bold C}$ to $\bmdr$.  
We say that a pair $(F, \zeta)$ 
is a {\it map of categories with $\V$-faces} 
from $({\bold C}, \Omega)$ to $({\bold C}^{\prime}, \Omega^{\prime})$ 
if $F\!: {\bold C} \to {\bold C}^{\prime}$ is a functor and  
$\zeta\!:\Omega \cong \Omega^{\prime} \circ F$ is a natural isomorphism.  
%
%
For categories with $\V$-faces
$({\bold C}_1, \Omega_1)$, $({\bold C}_2, \Omega_2)$, we denote
the category ${\bold C}_1 \times {\bold C}_2$ with $\V$-face
$\Omega_1 \otimes_R \Omega_2\!:$ $(X_1, X_2)$ $\mapsto $
$\Omega_1 (X_1) \otimes_R \Omega_2 (X_2)$ 
$(X_1 \in \ob {\bold C}_1, X_2 \in \ob {\bold C}_2)$ by 
$({\bold C}_1, \Omega_1) \times ({\bold C}_2, \Omega_2)$.

For a category with $\V$-face $({\bold C}, \Omega)$, 
we denote by $C ( \Omega ) = C({\bold C}, \Omega)$ 
the {\it coend} of $\Omega^* \otimes \Omega$, 
where $\Omega^* \otimes \Omega$ is viewed as a functor from 
${\bold C}^{\mathrm{op}} \times {\bold C}$ to $\bmde$ 
(cf. \cite{Maclane}).
By definition, these exists a unique family of maps 
$\kappa_X = \kappa_X^{\Omega}\!:$ 
$\Omega (X)^* \otimes \Omega (X) \to C(\Omega)$
$(X \in \mathrm{ob} {\bold C})$
in $\bold{Bmd} (E)$ such that (i) $\{\kappa_X \}_X$ 
is a dinatural transformation, that is, the diagram  
\begin{equation}
 \begin{CD}
 \Omega (Y)^* \otimes \Omega (X) @> \mathrm{id} \otimes \Omega (f)>> 
 \Omega (Y)^* \otimes \Omega (Y)\\
 @V\Omega (f)^* \otimes \mathrm{id}VV @VV\kappa_YV \\
 \Omega (X)^* \otimes \Omega (X)
 @>\kappa_{X}>> 
 C (\Omega)
 \end{CD}
\end{equation}
is commutative for each $f\!: X \to Y$
and that (ii) for another dinatural transformation
$\{ \iota_X\!:$ 
$\Omega (X)^* \otimes \Omega (X) \to M \}_X$,
there exists a unique map $h\!: C(\Omega) \to M$ 
such that $\iota_X = h \circ \kappa_X$ for each 
$X \in \mathrm{ob} {\bold C}$.
By the explicit description of the coend
via direct sums and a coequalizer, 
$C (\Omega)$ agrees with the coend of 
$\bar{\Omega}^* \otimes \bar{\Omega}\!:$
${\bold C}^{\mathrm{op}} \times {\bold C}$ $\to$ $\bold{Mod} (\K)$
as a vector space,
where $\bar{\Omega}$
denotes the composition of 
$\Omega$ with the forgetful functor 
$\bmdr \to \bold{mod} (\K)$.
Hence $C ( \Omega)$ becomes a coalgebra
(see e.g. Schauenburg \cite{SchauenburgTD} p. 29).


%
For $X \in \mathrm{ob} {\bold C}$, we define the map
$\delta_X = \delta_X^{\Omega}\!:$
$\Omega (X) \to \Omega (X) \otimes C (\Omega)$
by 
\begin{equation}
 \delta_X (m) = 
 \sum_i m_i \otimes \kappa_X (m^i \otimes m)
 \quad (m \in \Omega (X)),
\end{equation}
where $\{ m_i \}$ denotes a basis of $\Omega (X)$ and $\{ m^i \}$ denotes
its dual basis. 
Then $\delta_X$ gives a right $C(\Omega)$-comodule structure
on $\Omega (X)$, and satisfies 
\begin{gather}
 \delta_X (\lambda m \mu)  =
 \sum_{(m)} m_{(0)} \otimes \lambda m_{(1)} \mu,\\
 \sum_{(m)} \lambda m_{(0)} \mu \otimes m_{(1)} =
 \sum_{(m)} m_{(0)} \otimes \lamma m_{(1)} \muma
\end{gather}
%
%
%
for each $m \in \Omega (X)$ and $\lambda, \mu \in \V$.
Using these, we see that $C(\Omega)$ is a $\V$-dressed coalgebra. 
For example, if ${\bold C} = {\bold 1}$ has the only one object $I$ and 
the only one morphism, and if $\Omega = R$ sends $I$ to $R$, then, 
$C(\Omega)$ is isomorphic to $E$ as a $\V$-dressed coalgebra.

For categories with $\V$-faces $({\bold C}_1, \Omega_1)$
and $({\bold C}_2, \Omega_2)$, there exists a unique isomorphism 
$\phi_2\!: C({\bold C}_1, \Omega_1) \tensE C({\bold C}_2, \Omega_2)$
$\cong$
$C(({\bold C}_1, \Omega_1) \times ({\bold C}_2, \Omega_2))$
of $\V$-dressed coalgebras such that
\begin{equation*}
 \begin{CD}
 (\Omega_1 (X_1)^* \otimes \Omega_1 (X_1)) \otimes_E
 (\Omega_2 (X_2)^* \otimes \Omega_2 (X_2)) 
 @>\kappa^{\Omega_1}_{X_1} \otimes_E \kappa^{\Omega_2}_{X_2}>> 
 C (\Omega_1) \tensE C (\Omega_2) \\
 @V{\cong}VV @VV{\phi_2}V\\
 (\Omega_1 (X_1) \otimes_R \Omega_2 (X_2))^* \otimes
 (\Omega_1 (X_1) \otimes_R \Omega_2 (X_2))
 @>{(\kappa^{\Omega_1 \otimes_R \Omega_2})_{(X_1, X_2)}}>> 
 C(\Omega_1 \otimes_R \Omega_2)
 \end{CD}
\end{equation*}
is commutative for each
$X_1 \in \ob {\bold C}_1$ and $X_2 \in \ob {\bold C}_2$.

Suppose ${\bold C}$ is a $\K$-linear abelian category and $\Omega$ is 
a faithful $\K$-linear exact functor.
Then applying a primitive version of the Tannaka-Krein duality theorem to 
$\bar{\Omega}$, we see that $\Omega$ gives a categorical equivalence 
from ${\bold C}$ onto $\bold{com} (C(\Omega))$
(see e.g. Schauenburg \cite{SchauenburgTD} p. 39).

For a map 
$(F, \zeta)\!: ({\bold C}, \Omega) \to ({\bold C}^{\prime}, \Omega^{\prime})$ 
of categories with $\V$-faces, there exists a unique map
$C (F, \zeta)\!: C (\Omega) \to C (\Omega^{\prime})$ of
$\V$-dressed coalgebras such that the following diagram is commutative
for each $X \in \ob {\bold C}$: 
\begin{equation}
 \begin{CD}
 \Omega (X) @>\delta^{\Omega}_X>> \Omega (X) \otimes C (\Omega) \\
 @V\zeta_XVV @VV\zeta_X \otimes C(F, \zeta)V \\
 \Omega^{\prime} (F(X)) 
 @>\delta^{\Omega^{\prime}}_{F(X)}>> 
 \Omega^{\prime} (F(X)) \otimes C (\Omega^{\prime}).
 \end{CD}
\end{equation}
For another map 
$(F^{\prime}, \zeta^{\prime})\!:  ({\bold C}^{\prime}, \Omega^{\prime})$ $\to$
$({\bold C}^{\prime \prime}, \Omega^{\prime \prime})$
of categories with $\V$-face, we have 
$C((F^{\prime}, \zeta^{\prime}) \circ (F, \zeta)) =$
$C(F^{\prime}, \zeta^{\prime}) \circ C(F, \zeta)$,
where $(F^{\prime}, \zeta^{\prime}) \circ (F, \zeta)$ 
denotes the composition
$(F^{\prime} \circ F, \{ \zeta^{\prime}_{F(X)} \circ \zeta_X \}_X))$.

Let $({\bold C}, \Omega)$ be a category with $\V$-face such that 
${\bold C}$ $=$ $({\bold C}, \otimes)$ is a monoidal 
category with unit object $I$ and that
$\Omega$ is a monoidal functor with natural isomorphism 
$(\varphi_2)_{XY}\!: \Omega(X) \otimes_R \Omega(Y) \cong \Omega(X \otimes Y)$
and with isomorphism $\varphi_0\!: R \cong \Omega (I)$. 
Since $(\otimes, \varphi_2)\!:$
$({\bold C}, \Omega) \times ({\bold C}, \Omega) \to ({\bold C}, \Omega)$
is a map of categories with $\V$-faces, we obtain a map 
$m = C (\otimes, \varphi_2) \circ \phi_2\!:
C(\Omega) \otimes_E C(\Omega) \to C(\Omega)$ of $\V$-dressed coalgebras.
Similarly, considering $(I, \varphi_0)$ as a map of 
categories with $\V$-faces from $({\bold 1}, R)$ to $({\bold C}, \Omega)$, 
we obtain a map $\eta = C (I, \varphi_0)$ of $\V$-dressed coalgebras
from $E$ to $C(\Omega)$.
The triple $(C(\Omega), m, \eta)$ becomes a monoid in 
$\bold{DrCoalg_{\V}}$, 
and therefore it gives a $\V$-face algebra. 
Moreover, $\Omega$ gives a monoidal functor from ${\bold C}$ to 
$\bold{com}(C(\Omega))$.

Next, suppose that each object $X$ of ${\bold C}$
has a left dual $X^{\vee}$.
Using the uniqueness of the left dual of $\Omega (X)$ in $\bmdr$, 
we obtain a natural isomorphism 
$j_X\!: \Omega (X)^{\vee} \to \Omega (X^{\vee})$ 
in a natural manner.
Since
$({}^{\vee}, j)\!: ({\bold C}^{\op}, \Omega^{\vee})
\to ({\bold C}, \Omega)$ is a map of categories with $\V$-faces,
there exists a unique map
$S\!: C (\Omega)^{\bop} \to C (\Omega)$ of 
$\V$-dressed coalgebras such that
\begin{equation}
 \begin{CD}
 \Omega (X)^{\vee} @>{\mathrm{tw} \circ (\delta^{\Omega}_X)^{\flat}}>> 
 \Omega (X)^{\vee} \otimes C (\Omega)^{\bop} \\
 @V{j_X}VV @VV{j_X \otimes S}V \\
 \Omega (X^{\vee}) 
 @>\delta^{\Omega}_{X^{\vee}}>> 
 \Omega (X^{\vee}) \otimes C (\Omega)
 \end{CD}
\end{equation}
is commutative for each $X \in \ob {\bold C}$, 
where 
\begin{gather}
 (\delta^{\Omega}_X)^{\flat}(m^{\vee}) =
 \sum_i \sum_{(m_i)} \langle m^{\vee},\, {m_i}_{(0)} \rangle
 {m_i}_{(1)} \otimes m^i,\\
 \mathrm{tw} (c \otimes m^{\vee}) =
 m^{\vee} \otimes c
 \quad (m^{\vee} \in \Omega (X)^{\vee}, c \in C(\Omega)^{\bop}). 
\end{gather}
The map $S$ gives an antipode of $C(\Omega)$.

Let $({\bold C}, \otimes)$ be a monoidal category with unit object $I$.
We say that ${\bold C}$ is a {\it tensor category} over a field $\K$ if 
${\bold C}$ is a $\K$-linear abelian category and $\otimes$ is 
a bi-additive functor.
As the way stated above, the Tannaka-Krein duality is generalized as follows.

\begin{thm} 
 Let ${\bold C}$ be an essential small tensor category over $\K$
 and $\Omega\!: {\bold C} \to \bold{bmd}(R_{\V})$ 
 a faithful $\K$-linear exact monoidal functor.
 Then the coend $C (\Omega)$ becomes a $\V$-face algebra
 and $\Omega$ gives an equivalence 
 ${\bold C}$ $\cong$ $\bold{com}(C(\Omega))$
 of tensor categories.  
 If ,in addition, ${\bold C}$ is rigid, then $C (\Omega)$ has
 a bijective antipode. 						
\end{thm}

\noindent
{\it Remark}. $\!$
In \cite{fa}, we introduced the notion of an $R$-{\it face algebra}
for each commutative separable algebra $R$ over a field $\K$. 
Similarly to the``$\V$-Tannaka'' duality stated above,
we also have an $R$-Tannaka duality for arbitrary $R$.

\section{The canonical fiber functor} 
A tensor category ${\bold C}$ is called a {\it split finite semisimple category}
if there exists a finite set of objects $\{ L_{\lambda} | \lambda \in \V\}$
such that each object of ${\bold C}$ is isomorphic to a finite direct sum of
$L_{\lambda}$'s and that these objects satisfy the following Schur's
lemma;
\begin{equation}
 {\bold C} (L_{\lambda}, L_{\mu}) \cong 
 \K \delta_{\lambda \mu} \mathrm{id}_{L_{\lambda}}.
\label{SchurLem}
\end{equation}
Let ${\bold C}$ be a split finite semisimple tensor category.
For each object $X$ of ${\bold C}$, we define a bimodule 
$\Omega_0 (X)$ of $R = R_{\V}$ by setting
\begin{equation}
\label{Omega0Xdef}
 \lambda \Omega_0 (X) \mu = {\bold C} (L_{\mu}, L_{\lambda} \otimes X) 
 \quad (\lambda, \mu \in \V).
\end{equation}
Then $\Omega_0$ gives a faithful $\K$-linear exact functor 
from ${\bold C}$ to $\bmdr$ by setting
\begin{equation}
\label{Omega0fdef}
 \Omega_0 (f) (k) = ({\mathrm{id}}_{L_{\lambda}} \otimes f) \circ k 
\end{equation}
for each $f \in {\bold C} (X, Y)$ and 
$k \in  \lambda \Omega_0 (X) \mu$ $=$ 
${\bold C} (L_{\mu}, L_{\lambda} \otimes X)$. 
Moreover $\Omega_0$ gives a $\K$-linear exact faithful
monoidal functor by setting
\begin{equation*}
 (\varphi_2)_{X,Y}:\! \Omega_0 (X) \otimes_R \Omega_0 (Y) \to
 \Omega_0 (X \otimes Y) ;
\end{equation*}
\begin{gather}
\label{phi2def}
 f \otimes g \mapsto a_{L_{\lambda} X Y} \circ
 (f \otimes \mathrm{id}_Y) \circ g \nonumber\\ 
 (f \in \lambda \Omega_0 (X) \nu,\, 
 g \in \nu \Omega_0 (Y) \mu,\, \lambda, \mu, \nu \in \V),
\end{gather}
\begin{equation}
\label{phi0def}
 \varphi_0:\! R \to \Omega_0 (I); 
\quad  
 \lambda \mapsto r_{L_{\lambda}}^{-1}
 \quad (\lambda \in \V),
\end{equation}
where $a_{XYZ}:\! (X \otimes Y) \otimes Z$
$\to$ $X \otimes (Y \otimes Z)$ and 
$r_X:\! X \otimes I \to X$ denote the associativity
constraint and the right unit constraint of ${\bold C}$, respectively. 
Thus, applying the $\V$-Tannaka duality to $({\bold C}, \Omega_0)$, 
we obtain the following {\it canonical Tannaka duality}.

\begin{thm}
For each split finite semisimple category ${\bold C}$, there exists
a natural way to construct 
a $\V$-face algebra $C(\Omega_0)$ and an equivalence 
$\Omega_0\!:{\bold C} \cong \bold{com} (C(\Omega_0))$
of tensor categories.
\end{thm}
%
%
%
\noindent
{\it Example}. $\!$
Let $G$ be a finite group. 
%
%
%
The tensor category $\bold{com} (\C [G])$ is finite semisimple 
and $\{ \C g | g \in G \}$ gives a 
complete set of representatives of simple objects of $\bold{com} (\C [G])$.
The canonical Tannaka dual of this category has a basis
$\{ e^a_b [g] | a,b,g \in G \}$ which satisfies the following 
relations:
\begin{equation}
\label{eijpeklq} 
 e^a_b [g] e^c_d [h] = 
 \delta_{ag, c} \delta_{bg, d}
 e^a_b [gh],
\end{equation}
\begin{equation}
 \Delta( e^a_b [g] ) =
 \sum_{c \in G} e^a_c [g] \otimes e^c_b [g],
\quad
 \cu (e^a_b [g]) = \delta_{ab},
\end{equation}
\begin{equation}
 \stackrel{\scriptscriptstyle\circ}{a}\,
 =\, \sum_{b \in G} e^a_b [1], 
\quad
 b\, =\, \sum_{a \in G} e^a_b [1],
\end{equation}
\begin{equation}
 S (e^a_b [g]) = e^{bg}_{ag} [g^{-1}], 
\end{equation}
%
%
%
where $a,b,c,d,g,h \in G$.
%


\end{document}